\documentclass{scrartcl}
\usepackage[utf8]{inputenc}

\usepackage[english]{babel}
\usepackage{amsmath,amssymb}
\usepackage{enumitem} 
\usepackage{mathtools}
\usepackage{amsthm}
\usepackage{thmtools}
\usepackage{todonotes}
\usepackage{graphics}
\usepackage{multicol}
\usetikzlibrary{arrows,calc}

\usepackage{hyperref}

\declaretheorem[numberwithin=section, name=Lemma]{lemma}

\declaretheorem[name=Theorem,sibling=lemma]{theorem}
\declaretheorem[name=Corollary,sibling=lemma]{corollary}

\declaretheorem[style=definition, name=Example,sibling=lemma]{example}
\declaretheorem[style=definition,name=Definition,sibling=lemma]{definition}
\declaretheorem[style=definition, name=Remark,sibling=lemma]{remark}

\def\dotdiv{\mathbin{\ooalign{\hss\raise1ex\hbox{.}\hss\cr
  \mathsurround=0pt$-$}}}
\newcommand{\N}{\mathbb{N}}
\newcommand{\Z}{\mathbb{Z}}

\newcommand{\C}{{C}}

\newcommand{\T}{{T}}
\newcommand{\G}{{ G}}
\renewcommand{\H}{{ H}}
\renewcommand{\P}{{P}}
\newcommand{\AND}{\wedge}

\newcommand{\Pol}{\ensuremath{\operatorname{Pol}}}

\newcommand{\malt}{\Sigma_M}
\DeclareMathOperator{\Digraphs}{Digraphs}
\newcommand{\PosetDi}{\mathfrak P_{\Digraphs}}
\newcommand{\id}{{\operatorname{id}}}

\title{Maximal Digraphs With Respect to Primitive Positive Constructability\thanks{This research has been supported by the ERC, grant number xxx.}}
\author{Manuel Bodirsky and Florian Starke \\ Institut f\"ur Algebra, TU Dresden}
\date{March 2021}

\begin{document}

\maketitle

\begin{abstract}
    We study the class of all finite directed graphs (digraphs) up to primitive positive constructability. The resulting order has a unique greatest element, namely the digraph $P_1$ with one vertex and no edges. The digraph $P_1$ has a unique greatest lower bound, namely the digraph $P_2$ with two vertices and one directed edge. Our main result is a complete description of the greatest lower bounds of $P_2$; we call these digraphs \emph{submaximal}. We show that every digraph that is not equivalent to $P_1$ and $P_2$ is below one of the submaximal digraphs.
\end{abstract}

\section{Introduction}
A \emph{homomorphism} from a directed graph $G$ to a directed graph $H$ is a map from the vertices of $G$ to the vertices of $H$ which maps each edge of $G$ to an edge of $H$. Two directed graphs $G$ and $H$ are called \emph{homomorphically equivalent} if there is a homomorphism from $G$ to $H$ and from $H$ to $G$.
The study of the  \emph{homomorphism order}
on the class of all finite directed graphs (or short: \emph{digraphs}), factored by homomorphic equivalence, has a long history in graph theory. It is known to have a quite complicated structure; we refer to Ne\v{s}et\v{r}il and Tardif~\cite{NesetrilTardif} and the references therein. 

A classical topic in graph homomorphisms is the $H$-coloring problem, which is the computational problem of deciding whether a given finite digraph $G$ maps homomorphically to $H$. The computational complexity of this problem has been classified for finite undirected graphs $H$ by Hell and Ne\v{s}et\v{r}il~\cite{HellNesetril} in 1990: they are either in L or NP-complete. Feder and Vardi~\cite{FederVardi} proved that every finite-domain CSP is polynomial-time equivalent to an $H$-coloring problem for a finite \emph{directed} graph $H$\footnote{This result has been sharpened in~\cite{BulinDelicJacksonNiven}.}, and they conjectured
that each of these problems are either in P or NP-complete. 
This conjecture was eventually solved in 2017 by Bulatov and, independently, by Zhuk~\cite{BulatovFVConjecture,ZhukFVConjecture}. 
However, other long-standing open problems about the complexity of $H$-coloring for finite digraphs $H$ remain open, for example the characterisation of when this problem is in the complexity class L, or in NL~\cite{LinearDatalog,EgriLaroseTessonLogspace,Kazda18}. 

The border between polynomial-time tractable and NP-complete $H$-colouring problems can be described in terms of \emph{primitive positive (pp) constructions}, which is a concept that has been introduced by Barto, Opr\v{s}al, and Pinsker~\cite{wonderland} in the setting of general relational structures. The idea is that if $G$ has a pp construction in $H$, then, intuitively,  \emph{`$H$ can simulate $G$'}, and the $G$-coloring problem reduces (in logarithmic space) to the $H$-coloring problem. 
In particular, $H$-coloring is NP-hard if $K_3$ has a pp construction in $H$, where $K_3$ is the clique with three vertices, by reduction from the NP-hard three-colorability problem. It follows from the proof of the dichotomy conjecture that 
 otherwise $H$-coloring is in P. Note that pp constructability can also be used to
study the question of which $H$-coloring problems are in L or in NL. The surprising power of pp constructions is the motivation for studying pp constructions on finite digraphs more systematically. 

For digraphs $G$ and $H$ that have at least one edge, the definition of pp constructions takes the following elegant combinatorial form: 
$G$ pp constructs $H$ if there exists a digraph $K$  and $a,b \in V(K)^d$ for some $d \in {\mathbb N}$ 
 such that $G$ is homomorphically equivalent to the digraph with vertices $V(H)^d$ 
and where $(u,v)$ forms an edge if there is a homomorphism from $K$ to $H$ that maps $(a,b)$ to $(u,v)$. 
We write $H \leq G$ if $G$ has a pp construction in $H$. It can be shown that $\leq$ is transitive (Corollary 3.10 in~\cite{wonderland}) and so it gives rise to a partial order $\PosetDi$ on the class of all finite digraphs (where we take the liberty to identify two digraphs $G$ and $H$ if they pp construct each other). 
Since all finite digraphs have a pp construction in $K_3$ (see, e.g. [11]), it is the smallest element of the poset $\PosetDi$.
The poset also has a greatest element, namely the digraph
$P_1$ with just one vertex and no edges. 
The digraph $P_1$ has a unique greatest lower bound in $\PosetDi$, namely the digraph $P_2$ consisting of two vertices and one directed edge; this is not hard to see and will be shown in Section~\ref{sect:result}. 

\begin{figure}
    \centering
    \begin{tikzpicture}[scale=1.3]
    \node (0) at (2,2)  {$\P_{1} \equiv C_{1}$};
    \node (1) at (2,1)  {$\P_{2}$};
    \node (20) at (0,0)  {$\T_{3}$};
    \node (21) at (1,0)  {$\C_{2}$};
    \node (22) at (2,0)  {$\C_{3}$};
    \node (23) at (3,0)  {$\C_{5}$};
    \node (24) at (4,0)  {$\dots$};

    \node (3) at (2,-0.7)  {$\vdots$};
    \node[rotate = 45] (3) at (0.8,-0.7)  {$\vdots$};
    \node[rotate = -45] (3) at (3.2,-0.7)  {$\vdots$};

    \node (4) at (2,-1.5)  {$K_3$};
    
    \path
        (0) edge (1)
        (1) edge (20)
        (1) edge (21)
        (1) edge (22)
        (1) edge (23)
        (1) edge (24)
        ;
    
    \end{tikzpicture}
    \caption{The pp constructability poset on finite digraphs.}
    \label{fig:main}
\end{figure}
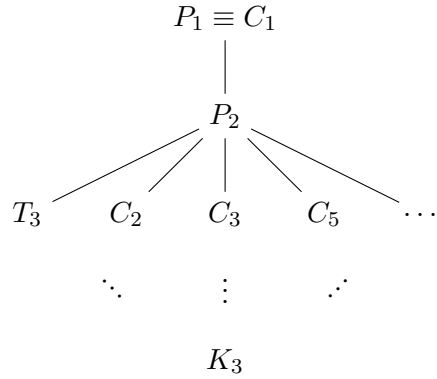

In this article, we present a complete description of the greatest lower bounds of $P_2$ in $\PosetDi$; we call these digraphs \emph{submaximal}. We also prove that every finite digraph which does not pp constructs $P_2$ is smaller than one of the submaximal digraphs (Theorem~\ref{thm:submaximalGraphs}; also see Figure~\ref{fig:main}). 
The submaximal digraphs are:
\begin{itemize}
    \item The directed cycles $C_p$ for $p$ prime.
    (For $k \in {\mathbb N}^+$, the directed cycle $C_k$ 
    is defined to be the digraph 
    $({\mathbb Z}_{k};\{(u,v) \mid u+1=v \mod k\})$.) 
    \item $T_3 \coloneqq (\{0,1,2\},<)$, the transitive tournament with three vertices.
\end{itemize}

\subsection*{Related work}
The pp constructability poset for smooth digraphs, i.e., digraphs where every vertex has indegree at least one and outdegree at least one (digraphs without sources and sinks), has been described in~\cite{smooth-digraphs}. 
The pp constructability poset on general relational structures over a two-element set has been described  in~\cite{PPPoset}.


\section{Minor conditions} 
Primitive positive constructability has a universal algebraic characterisation; this characterisation plays a role in our proof, so we present it here. 
If $H = (V,E)$ is a digraph then $H^k$ denotes the $k$-th direct power of $H$, which is the digraph with vertex set $V^k$ and edges set $$\{((u_1,\dots,u_k),(v_1,\dots,v_k)) \mid (u_1,v_1) \in E,\dots,(u_k,v_k) \in E\}.$$
A \emph{polymorphism} of $H$ is a homomorphism $f$ from $H^k$ to $H$, for some $k \in {\mathbb N}$, which is called the \emph{arity} of $f$. We write $\Pol(H)$ for the set of all polymorphisms of $H$. This set contains the projections and is closed under composition.\footnote{Sets of operations with these properties are called \emph{clones} in universal algebra.} An operation $f$ is called \emph{idempotent} if $f(x,\dots,x) = x$ for all $x \in V$. 

A central topic in universal algebra are \emph{minor conditions}. If $f \colon V^k \to V$ is an operation and $\sigma \colon \{1,\dots,k\} \to \{1,\dots,n\}$ is a function,
then $f_\sigma$ denotes the operation 
$$(x_1,\dots,x_n) \mapsto f(x_{\sigma(1)},\dots,x_{\sigma(k)}),$$
and $f_\sigma$ is called a \emph{minor} of $f$. 
A \emph{minor condition} is a set $\Sigma$ of expressions of the form $f_{\sigma} = g_{\tau}$ where $f$ and $g$ are function symbols ($f$ and $g$ might be the same symbol).  

\begin{example}
An operation $f \colon V^n \to V$ is called \emph{cyclic}
if for all $x_1,\dots,x_n \in V$
$$f(x_1,x_2,\dots,x_n) = f(x_2,\dots,x_n,x_1).$$ This condition 
can be expressed by the minor condition 
$$\Sigma_n \coloneqq \{f_{\id} = f_{\tau}\}$$ where $\id$ denotes the identity function on $\{1,2,\dots,n\}$ and $\tau$ denotes the cyclic permutation $(1,2,\dots,n)$ on $\{1,\dots,n\}$. 
\end{example} 

If a minor condition $\Sigma$ contains several expressions, then 
different expressions in $\Sigma$ might share the same function symbols.

\begin{example}
An idempotent operation $f$ is called a \emph{Maltsev operation} 
if for all $x,y \in V$ 
\[f(y,y,x) = f(x,x,x) = f(x,y,y) .\] 
This condition can be expressed by the minor condition 
\[\Sigma_M \coloneqq \{f_{\sigma} = f_\tau, f_{\tau} = f_{\rho}\}\]
where $\sigma, \tau, \rho \colon \{1,2,3\} \to \{1,2\}$ are given by $\sigma(1,2,3) = (2,2,1)$,  $\tau(1,2,3) = (1,1,1)$, and 
$\rho(1,2,3) = (1,2,2)$. 
\end{example}

A set of operations $F$ \emph{satisfies} a minor condition $\Sigma$ if the function symbols in $\Sigma$ can be replaced by operations from $F$ so that all the expressions in $\Sigma$ hold; in this case we write $F \models \Sigma$. 
If $H$ is a digraph, then $\Sigma(H)$ denotes the class of all minor conditions that are satisfied in $\Pol(H)$. 

\begin{theorem}[Barto, Opr\v{s}al, and Pinsker~\cite{wonderland}]\label{thm:wonderland} 
Let $G$ and $H$ be finite digraphs. Then 
\begin{align*}
    H \text{ pp constructs } G &&\text{if and only if} && \Sigma(H)\subseteq \Sigma(G). 
\end{align*}
\end{theorem}

\section{The pp construction poset}\label{sect:result}
We have already defined pp constructability for digraphs in the introduction, but present an equivalent description here which is convenient when specifying pp constructions, and which is also closer to the presentation of Barto, Opr\v{s}al, and Pinsker~\cite{wonderland}. 
A \emph{primitive positive formula} is a formula $\phi(x_1,\dots,x_k)$ of the form
$$ \exists y_1,\dots,y_n (\psi_1 \wedge \cdots \wedge \psi_m)$$
where each of the formulas $\psi_1,\dots,\psi_m$ is of the form $\bot$ (for \emph{false}), of the form $z_1=z_2$, or of the form $E(z_1,z_2)$ where
$z_1,z_2$ are variables from $\{x_1,\dots,x_k,y_1,\dots,y_n\}$.

\begin{definition}
Let $H = (V,E)$ be a digraph. A digraph $G$ with vertex set $V^d$ is called a \emph{pp power of $H$ of dimension $d$} if there exists a primitive positive formula $\phi(x_1,\dots,x_d,y_1,\dots,y_d)$
such that 
the edge set of $G$ equals 
$$\{((u_1,\dots,u_d),(v_1,\dots,v_d)) \mid \phi(u_1,\dots,u_d,v_1,\dots,v_d) \text{ holds in } H\}.$$
\end{definition}
It follows from the definitions that $H \leq G$ if and only if $G$ is homomorphically equivalent to a pp power of $H$. 
We write 
\begin{itemize}
\item $H \equiv G$ if $H \leq G$ and $G \leq H$;
\item 
$H < G$ if $H \leq G$ and not $G \leq H$. 
\end{itemize}

\begin{lemma}
$\P_1$ is the greatest element of $\PosetDi$. Moreover, $\P_1 \equiv \C_1$. 
\end{lemma}
\begin{proof}
Let $\G$ be a finite digraph. Consider the pp power of $\G$ of dimension one given by the formula $\phi(x,y) \coloneqq \;\perp$. The resulting digraph has no edges and is therefore homomorphically equivalent to $\P_1$.
Now consider the pp power of $\G$ of dimension one given by the formula $\phi(x,y) \coloneqq (x=y)$. 
The resulting digraph is homomorphically equivalent to the digraph $\C_1$ with one vertex and a loop, which implies the statement. 
\end{proof}

In the proof of the following lemma we need the fundamental concept of \emph{cores} from the theory of graph homomorphisms (see, e.g.,~\cite{HNBook}). A digraph $H = (V,E)$ is called a \emph{core} if every \emph{endomorphism} of $H$ (i.e., every homomorphism from $H$ to $H$) is an \emph{embedding} (i.e., an
isomorphism between $H$ and an induced subgraph of $H$). It is easy to see that every finite  digraph $H$ is homomorphically equivalent to a core digraph, and that all core digraphs $G$ that are homomorphically equivalent to $H$ are isomorphic to each other; we therefore call $G$ \emph{the} core of $H$. 
When studying $\PosetDi$ we may therefore restrict our attention to core digraphs; the big advantage of cores is the following useful lemma.

\begin{lemma}[follows from Lemma 3.9 in~\cite{wonderland}]\label{lem:constants}
Let $H = (V,E)$ be a finite core digraph. Then $H \leq G$ if and only if $G$ is homomorphically equivalent to a pp power of $H$ where the primitive positive formula might additionally contain conjuncts of the form $x = c$ where $x$ is a variable and $c \in V$ is a constant. 
\end{lemma}

\begin{lemma}
We have $\P_2 < \P_1$. Moreover, 
$\P_2$ is the only coatom of $\PosetDi$, 
i.e., $\P_2$ is the unique greatest lower bound of $\P_1$ in $\PosetDi$. 
\end{lemma}
\begin{proof}
We have already seen that $\P_2 \leq \P_1$. 
To prove that $\P_2 \not \leq P_1$, 
first observe that $\P_1$ has constant  polymorphisms, while $\P_2$ does not. Let $\Sigma_c \coloneqq \{f_{\rho} = f_{\sigma} \}$ 
where $f$ is a unary function symbol, 
$\rho \colon \{1\} \to \{1,2\}$ is constant 1, and $\sigma \colon \{1\} \to \{1,2\}$ is constant 2.  
Then $\Pol(\P_1) \models \Sigma_c$, 
but $\Pol(\P_2) \not \models \Sigma_c$. Then (the easy direction of) 
Theorem~\ref{thm:wonderland} implies that $\P_1 \leq \P_2$ does not hold. 

For the second statement, let $\G$ be a finite digraph such that $\G<\P_1$. 
We have to show that $\G \leq \P_2$. 
Without loss of generality we may assume that $\G$ is a core. Hence, by Lemma~\ref{lem:constants}, we can use constants in pp constructions. Note that $G$ must have at least two different vertices $u$ and $v$. The pp power of $\G$ of dimension one given by the formula $\phi(x,y) \coloneqq (x=u) \AND (y=v)$ is a digraph that has exactly one edge that is not a loop and that is therefore homomorphically equivalent to $\P_2$.
\end{proof}

The following theorem is our main result and will be shown in the remainder of the article; see Figure~\ref{fig:main}. 

\begin{theorem}\label{thm:submaximalGraphs}
The submaximal elements of $\PosetDi$ are precisely $\T_3$, $\C_2$, $\C_3$, $\C_5$, $\dots$ 
If $\G$ is a finite digraph that does not have a pp construction in $P_2$, then $\G\leq \T_3$ or $\G\leq\C_p$ for some prime $p$.
\end{theorem}

\section{Submaximal digraphs and minor conditions}
We first discuss which of 
the minor conditions that we have encountered 
are satisfied by the polymorphisms of 
the digraphs that appear in Theorem~\ref{thm:submaximalGraphs}. 
The following facts are well-known; we present the proof for the convenience of the reader. 

\begin{lemma}\label{lem:CpConditions1}
Let $p$ and $q$ be primes. Then
$\Pol(C_p) \models \Sigma_q$ if and only if $q \neq p$. 
\end{lemma}
\begin{proof}
 If $p\neq q$, then there is an $n\in \N^+\!$ such that $p\cdot n=1\pmod q$. The map
\[(x_1,\dots,x_p) \mapsto n\cdot(x_1 +\ldots + x_p) \pmod q\]
is a polymorphism of $\C_q$ satisfying $\Sigma_p$.

Now suppose that $p=q$. We assume for contradiction that $f$ is a polymorphism of $\C_q$ satisfying $\Sigma_p$. Then
\[f(0,\dots,q-2,q-1) = a = f(1,\dots,q-1,0)\]
and hence $(a,a) \in E$, which is impossible since $C_p$ has a loop only if $p=1$. 
\end{proof}

\begin{lemma}\label{lem:CpConditions2}
$\Pol(C_n) \models \Sigma_M$ for every $n \in {\mathbb N}$. 
\end{lemma}
\begin{proof}
The ternary operation $(x_1,x_2,x_3) \mapsto x_1 - x_2 + x_3 \pmod n$ is a Maltsev polymorphism of $C_n$. 
\end{proof}

\begin{figure}
    \centering
    \begin{tikzpicture}[scale=1.3]
    \node (a) at (0,0) {$a$};
    \node (b) at (0,1) {$b$};
    \node (c) at (1,1) {$c$};
    \node (d) at (1,0) {$d$};
    \path[->]
        (a) edge (b)
        (c) edge (b)
        (c) edge (d)
        (a) edge[dashed] (d)
        ;
    \end{tikzpicture}
    \caption{Rectangularity in digraphs.}
    \label{fig:rect}
\end{figure}
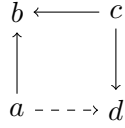

A finite digraph $H$ is called 
\emph{$k$-rectangular} if whenever $H$ contains directed paths of length $k$ from $a$ to $b$, from $c$ to $b$, 
and from $c$ to $d$, then also from $a$ to $d$. See Figure~\ref{fig:rect}. 
A digraph $H$ is called \emph{totally rectangular} if it is $k$-rectangular for all $k \geq 1$. 

\begin{lemma}\label{lem:maltIffRect}
A finite digraph $H$ is totally rectangular if and only if
it has a Maltsev polymorphism. 
A finite core digraph $H$ has a Maltsev polymorphism if and only if 
$\Pol(H) \models \malt$. 
\end{lemma}
\begin{proof}
The first part of the statement is Corollary 4.11 in~\cite{CarvalhoEgriJacksonNiven}. For the second statement, 
let $H = (V,E)$ be a core digraph which has a polymorphism $f$ that satisfies $f(x,y,y)=f(x,x,x)=f(y,y,x)$ for all $x,y \in V$; we have to find a polymorphism that is additionally idempotent. 
Note that the function $x \mapsto f(x,x,x)$ is an endomorphism; since $H$ is a core, the endomorphism is injective. Since $H$ is finite the endomorphism must in fact be an automorphism, and has an inverse $i$ which is an endomorphism as well. Then the operation $(x_1,x_2,x_3) \mapsto i(f(x_1,x_2,x_3))$ is idempotent and a Maltsev operation. 
\end{proof}

\begin{lemma}\label{lem:T3Conditions}
$\Pol(T_3) \models \Sigma_n$ for every $n \in {\mathbb N}$, but $\Pol(T_3) \not \models \Sigma_M$. 
\end{lemma}
\begin{proof}
The operation $(x_1,\dots,x_n) \mapsto \max(x_1,\dots,x_n)$ is a polymorphism of $T_3$ that satisfies $\Sigma_n$. 
On the other hand, $T_3 = (\{0,1,2\},E)$ is not 1-rectangular, witnessed by $(1,2),(0,2),(0,1) \in E$ but $(1,1) \notin E$; the second  statement therefore follows from 
Lemma~\ref{lem:maltIffRect}. 
\end{proof}


The following theorem states that 
the digraph $\P_2$ is the unique smallest element of 
$\PosetDi$ that satisfies $\Sigma_M$ and $\Sigma_p$ for all $p$ prime. 

\begin{theorem}\label{thm:maltAndCyclImplyIdempotent}
Let $\G$ be a finite digraph
that satisfies $\malt$
and $\Sigma_p$ for all primes $p$. Then 
$\P_2 \leq \G$. 
\end{theorem}

In the proof of Theorem~\ref{thm:maltAndCyclImplyIdempotent} we make use the following result of Carvalho, Egri, Jackson, and Niven~\cite{CarvalhoEgriJacksonNiven}, which guides us in our further proof steps. 

\begin{theorem}[Lemma 3.10 in~\cite{CarvalhoEgriJacksonNiven}]\label{thm:maltImpliesPathOrDuoc}
If $\G$ is totally rectangular, then $\G$ is homomorphically equivalent to either a directed path or a disjoint union of directed cycles.
\end{theorem}

We write $\P_k$ for the directed path with the $k$ vertices $\{0,\dots,k-1\}$. 

\begin{lemma}\label{lem:IdempConstructPaths}
The digraph $\P_2$ pp constructs $\P_{k}$ for all $k\in\N^+$\!.
\end{lemma}
\begin{proof}
Clearly, $\P_2 \leq \P_1$ and $\P_2 \leq \P_2$.
Let $k\geq3$ and consider the pp power $\G$ of $\P_2$ of dimension $k-1$ given by the following formula 
$\phi(x_1,\dots x_{k-1},y_1,\dots,y_{k-1})$ 
\begin{align*}
    (x_1 = y_2) \AND (x_2=y_3) \AND \dots \AND (x_{k-2}=y_{k-1}) \AND E(x_{k-1},y_1). 
\end{align*}
Then $\G$ contains the following path of length $k$
\[(0,0,\dots,0)\to (1,0,\dots,0) \to (1,1,\dots,0) \to \dots \to (1,1,\dots,1).\]
which shows that there exists a homomorphism from $\P_k$ to $\G$. Note that whenever there is an edge from $u$ to $v$ in $\G$, then the tuple $v$ contains exactly one 1 more than the tuple $u$. 
Therefore, the function $V(\G) \to \{1,\dots,k\}$ that maps $v$ to the number of 1's in $v$ is a homomorphism from $\G$ to $\P_k$. Hence $\P_2\leq\P_k$.
\end{proof}

\begin{proof}[Proof of Theorem~\ref{thm:maltAndCyclImplyIdempotent}]
Let $\G$ be a finite digraph satisfying $\malt$ and $\Sigma_p$ for every prime $p$. By Lemma~\ref{lem:maltIffRect} and Theorem~\ref{thm:maltImpliesPathOrDuoc} there are two cases to consider:
the first is that $\G$ is homomorphically equivalent to $\P_k$ for some $k$. Then $\P_2\leq \G$ by Lemma~\ref{lem:IdempConstructPaths}. 

The second case is that $\G$ is homomorphically equivalent to  a disjoint union of directed cycles.
Without loss of generality we may assume that $\G$ is a disjoint union of directed cycles.   
Let $(a_0,\dots,a_{\ell-1})$ be a shortest cycle in $\G$. Let $p$ be a prime and $k\in\N^+\!$ such that $p\cdot k=\ell$, and
let $f\in\Pol(\G)$ be a function that witnesses that  $\Pol(\G)\models\Sigma_p$. 
Then
\begin{align*}
    f(a_0,a_k,\dots,a_{(p-1)\cdot k})=a=f(a_k,a_{2 k},\dots,a_{0}).
\end{align*}
Since $f$ is a polymorphism of $\G$ there is a directed path of length $k$ from $a$ to $a$. Thus, $\G$ contains a directed cycle whose length divides $k$, which  contradicts the assumption that $\ell$ is the length of the shortest directed cycle in $\G$. Therefore, $\ell$ has no prime divisors, and $\ell=1$. So $\G$ contains a loop and hence is homomorphically equivalent to $C_1$; it follows that $\P_2 \leq \G$. 
\end{proof}

\section{Proof of the main result}

We use the following general result about when a finite digraph can pp construct a finite disjoint union of cycles.
\begin{lemma}[Lemma 6.8 in~\cite{smooth-digraphs}]
Let $\C$ be a finite disjoint union of cycles and let $\G$ be a finite digraph. Then 
\begin{align*}
\G\leq \C && \text{iff} && \Pol(\G)\models\Sigma_{C\dotdiv c} \text{ implies } \Pol(\C)\models\Sigma_{C\dotdiv c}\text{ for all $c\in\N^+$\!.} 
\end{align*}
\end{lemma}

For the special case that  $C=C_p$, there are only two conditions of the form $\Sigma_{C\dotdiv c}$, namely $\Sigma_1$, which is trivial, and $\Sigma_p$, which is not satisfied by $C_p$. Hence, we obtain the following result. 
\begin{theorem}
\label{thm:notSatSigma}
If $\G$ is a finite digraph. If $p$ is a prime number such that $\Pol(\G) \not\models\Sigma_p$, then $\G \leq\C_p$.
\end{theorem}


We also need a similar results for $\malt$ instead of $\Sigma_p$. 

\begin{lemma}\label{lem:notSatMalt}
Let $\G$ be a finite digraph. 
If $\Pol(\G) \not \models \malt$, then $\G \leq \T_3$. 
\end{lemma}

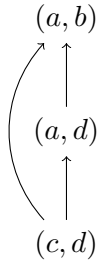
\begin{figure}
    \centering
    \begin{tikzpicture}
    \node (0) at (0,0) {$(c,d)$};
    \node (1) at (0,1.5) {$(a,d)$};
    \node (2) at (0,3) {$(a,b)$};
    
    \path[->]
        (0) edge (1)
        (1) edge (2)
        (0) edge[bend left=42] (2)
        ;
    \end{tikzpicture}
    \caption{The primitive positive construction of $T_3$ in the proof of Lemma~\ref{lem:notSatMalt}.}
    \label{fig:T3constr}
\end{figure}

\begin{proof}
Since $\leq$ is transitive we may assume without loss of generality that $\H = (V,E)$ is a core. By Lemma~\ref{lem:maltIffRect}, $\H$ is not totally rectangular. Hence, there are vertices $a,b,c,d \in V$ such that in $\G$ there are directed paths of length $k$ from $a$ to $b$, from $c$ to $b$, from $c$ to $d$, and there is no directed path of length $k$ from $a$ to $d$. Note that by Lemma~\ref{lem:constants} we are allowed to use constants in pp constructions. 
We write $x\stackrel{k}\to y$
as a shortcut for the primitive positive formula 
$\exists u_1,\dots,u_{k-2} (E(x,u_1) \wedge E(u_1,u_2) \wedge \cdots \wedge E(u_{k-2},y))$. 
Consider the pp power of $\G$ of dimension two given by the formula
\begin{align*}
    \phi(x_1,x_2,y_1,y_2)&\coloneqq x_1\stackrel{k}\to y_2 \AND (x_2=d) \AND (y_1=a).
\end{align*}

Let $\H$ be the resulting digraph. Consider the vertices $v_0=(c,d)$, $v_1=(a,d)$, and $v_2=(a,b)$ of $\H$. Note that the only vertex of $\H$ that can have incoming and outgoing edges is $v_1$. Since there is no path of length $k$ from $a$ to $d$ the vertex $v_1$ does not have a loop. Furthermore, $\H$ has the edges $(v_0,v_1), (v_1,v_2),$ and $(v_0,v_2)$ (see Figure~\ref{fig:T3constr}). 
Hence, $i\mapsto v_i$ is an embedding of $\T_3$ into $\H$. 
Let $V_0$ be the set of all vertices in $H$ that have outgoing edges and $V_2$ be the set of all vertices in $H$ that have incoming edges. Note that $V_0\mathbin\Delta V_2$ consists of $v_1$ and all isolated vertices. Clearly, $V_0\setminus V_2$, $V_0\mathbin\Delta V_2$, and $V_2\setminus V_0$ form a partition of $V(H)$ and
the map
\begin{align*}
    v\mapsto
    \begin{cases}
    v_2&\text{if }v\in V_2\setminus V_0\\
    v_1&\text{if }v\in V_0\mathbin\Delta V_2\\
    v_0&\text{if }v\in V_0\setminus V_2
    \end{cases}
\end{align*}
is a homomorphism from $\H$ to $\T_3$. Hence $\G \leq \T_3$.
\end{proof}

\begin{proof}[Proof of Theorem~\ref{thm:submaximalGraphs}]
Let $\G$ be a digraph such that $\P_2\not\leq \G$. 
Theorem~\ref{thm:maltAndCyclImplyIdempotent} implies that either $\Pol(\G)$ does not satisfy $\malt$ or that it does not satisfy $\Sigma_p$ for some prime $p$. In the first case $\G \leq \T_3$, by Lemma~\ref{lem:notSatMalt}. 
In the second case $\G\leq\C_p$, by Theorem~\ref{thm:notSatSigma}.
Hence, all submaximal elements of $\PosetDi$ are contained in 
$\{\T_3, \C_2 ,\C_3, \C_5, \dots\}$. 
Lemma~\ref{lem:CpConditions1}, Lemma~\ref{lem:CpConditions2},  and Lemma~\ref{lem:T3Conditions} imply that these digraphs form an antichain in $\PosetDi$, and hence 
each of these digraphs is submaximal. 
\end{proof}


Note that our result implies the following. 

\begin{corollary}\label{cor:maltAndCyclImplyIdempotent}
If a finite digraph $\G$ satisfies $\malt$, $\Sigma_2$, $\Sigma_3$, $\Sigma_5$, $\dots$, then any minor condition satisfied by $\Pol(\P_2)$ is also satisfied by $\Pol(\G)$.
\end{corollary}

The statement of Corollary~\ref{cor:maltAndCyclImplyIdempotent}  may also be phrased as
$$\{\malt, \Sigma_2, \Sigma_3, \Sigma_5, \dots\}\subseteq \Sigma(\G) \quad \Rightarrow \quad  \Sigma(\P_2)\subseteq \Sigma(\G).$$

\begin{remark}
We do not know whether Corollary~\ref{cor:maltAndCyclImplyIdempotent} holds for arbitrary clones of operations on a finite set, instead of just clones of the form $\Pol(\G)$ for a finite digraph $G$. However, the statement is false for clones of operations on an infinite set, as illustrated by the clone of operations on ${\mathbb Q}$ of the form
$(x_1,\dots,x_n) \mapsto a_1x_1+\cdots+a_nx_n$ for 
$a_1,\dots,a_n \in {\mathbb Q}$ such that $a_1 + \cdots + a_n = 1$. 
This clone satisfies $\Sigma_n$ for every $n \in \mathbb N$,
and contains the function $(x_1,x_2,x_3) \mapsto x_1 - x_2 + x_3$, so it also satisfies $\Sigma_M$. 
However, it is easy to see that it does not contain an operation $f$ that satisfies
$$f(x,x,y) = f(y,y,x) = f(x,y,y) = f(y,x,x)$$
for all $x,y \in {\mathbb Q}$;
however, this minor condition is satisfied by $\Pol(\P_2)$ (for example by $f = \max$). 
\end{remark}

\begin{remark}
Many, but not all the statements that we have shown also apply to \emph{infinite} digraphs. In Theorem~\ref{thm:wonderland}, only the forward direction holds if $G$ and $H$ are infinite; however, in this text we only used the forward direction of this theorem. 

Every digraph with a Maltsev polymorphism is totally rectangular even if the digraph
is infinite. 
The proof of  Theorem~\ref{thm:maltImpliesPathOrDuoc} of Carvalho, Egri, Jackson, and Niven can be generalised to show that 
every infinite digraph which is totally rectangular 
homomorphically equivalent to an infinite disjoint unions of cycles or one of the infinite paths $P^\infty\coloneqq(\N,\{(u,u+1)\mid u\in\N\})$, 
$P_\infty\coloneqq(\N,\{(u+1,u)\mid u\in\N\})$, 
the disjoint union $P_\infty + P^\infty$ of $P_\infty$ and $P^\infty$, 
and  $P^\infty_\infty\coloneqq(\Z,\{(u,u+1)\mid u\in\Z\})$. All of these graphs have a Maltsev polymorphism. 

Infinite disjoint unions of cycles are clearly not submaximal. 
Clearly, $P_2$ cannot pp construct the core digraphs $P_{\infty}$, $P^{\infty}$, $P_\infty + P^\infty$, and $P^\infty_\infty$, 
and these graphs can pp construct $P_2$. Clearly $P_{\infty}$ and $P^{\infty}$ pp-construct each other. 
We do not know whether these graphs are submaximal in the class of all digraphs.  
\end{remark}

\section{Concluding remarks}
Primitive positive constructability orders finite digraphs $H$ by their `strength' with respect to the $H$-coloring problem. Many deep combinatorial statements about graphs and digraphs can be phrased in terms of this order. We showed that at least the top region of the resulting poset can be described completely.
A full description of the entire poset $\PosetDi$ would be highly desirable.
We state three concrete open problems.
\begin{enumerate}
    \item Is $\PosetDi$ a lattice? (Primitive positive constructability is known to form a join meet-lattice on the class of all finite relational structures factored by homomorphic equivalence, but it is not clear to the authors whether the clone product construction for the meet used there can be carried out in the category of digraphs.) 
    \item Does $\PosetDi$ contain infinite ascending chains? (We have seen an infinite antichain in this article; an infinite descending chain of digraphs with a Maltsev polymorphism can be found in~\cite{smooth-digraphs} and 
    the existence of infinite descending chains of digraphs without a Maltsev polymorphism follows from results of~\cite{BMOOPW}, and also from results in~\cite{JacksonKN16}.) 
    \item What are the greatest lower bounds of $T_3$ in $\PosetDi$? 
\end{enumerate}

We already mentioned that the pp constructability poset on disjoint unions of cycles has been described in~\cite{smooth-digraphs}; in particular, it contains no infinite ascending chains and is a lattice. 
Note that this result combined with the result of the present paper shows that when exploring $\PosetDi$ it only remains to explore the interval between $K_3$ and $T_3$: 
if a digraph $\H$ does not have a Maltsev polymorphism, then we proved that it is below $T_3$ (and above $K_3$);
otherwise, it is homomorphically equivalent to a directed path or a disjoint union of cycles and hence falls into the region that has already been completely described.

\bibliographystyle{unsrt}
\bibliography{global}

\def\cprime{$'$} \def\cprime{$'$} \def\cprime{$'$}
\begin{thebibliography}{10}

\bibitem{NesetrilTardif}
J.~Ne\v{s}et\v{r}il and C.~Tardif.
\newblock Duality theorems for finite structures.
\newblock {\em Journal of Combininatorial Theory, Series B}, 80:80--97, 2000.

\bibitem{HellNesetril}
Pavol Hell and Jaroslav Ne\v{s}et\v{r}il.
\newblock On the complexity of {H}-coloring.
\newblock {\em Journal of Combinatorial Theory, Series B}, 48:92--110, 1990.

\bibitem{FederVardi}
Tom\'as Feder and Moshe~Y. Vardi.
\newblock The computational structure of monotone monadic {SNP} and constraint
  satisfaction: {a} study through {D}atalog and group theory.
\newblock {\em {SIAM} Journal on Computing}, 28:57--104, 1999.

\bibitem{BulinDelicJacksonNiven}
Jakub Bulin, Dejan Delic, Marcel Jackson, and Todd Niven.
\newblock A finer reduction of constraint problems to digraphs.
\newblock {\em Log. Methods Comput. Sci.}, 11(4), 2015.

\bibitem{BulatovFVConjecture}
Andrei~A. Bulatov.
\newblock A dichotomy theorem for nonuniform {CSP}s.
\newblock In {\em 58th {IEEE} Annual Symposium on Foundations of Computer
  Science, {FOCS} 2017, {B}erkeley, {CA}, {USA}, {O}ctober 15-17}, pages
  319--330, 2017.

\bibitem{ZhukFVConjecture}
Dmitriy~N. Zhuk.
\newblock A proof of {CSP} dichotomy conjecture.
\newblock In {\em 58th {IEEE} Annual Symposium on Foundations of Computer
  Science, {FOCS} 2017, {B}erkeley, {CA}, {USA}, {O}ctober 15-17}, pages
  331--342, 2017.
\newblock https://arxiv.org/abs/1704.01914.

\bibitem{LinearDatalog}
V\'{\i}ctor Dalmau.
\newblock Linear datalog and bounded path duality of relational structures.
\newblock {\em Logical Methods in Computer Science}, 1(1), 2005.

\bibitem{EgriLaroseTessonLogspace}
L{\'a}szl{\'o} Egri, Benoit Larose, and Pascal Tesson.
\newblock Symmetric datalog and constraint satisfaction problems in logspace.
\newblock In {\em Proceedings of the Symposium on Logic in Computer Science
  ({LICS})}, pages 193--202, 2007.

\bibitem{Kazda18}
Alexandr Kazda.
\newblock nnn-permutability and linear datalog implies symmetric datalog.
\newblock {\em Log. Methods Comput. Sci.}, 14(2), 2018.

\bibitem{wonderland}
Libor Barto, Jakub Opr\v{s}al, and Michael Pinsker.
\newblock The wonderland of reflections.
\newblock {\em Israel Journal of Mathematics}, 223(1):363--398, 2018.

\bibitem{smooth-digraphs}
Manuel Bodirsky, Florian Starke, and Albert Vucaj.
\newblock Smooth digraphs modulo primitive positive constructability.
\newblock {\em International Journal on Algebra and Computation (to appear)},
  2021.
\newblock Preprint available at ArXiv:1906.05699.

\bibitem{PPPoset}
Manuel Bodirsky and Albert Vucaj.
\newblock Two-element structures modulo primitive positive constructability.
\newblock {\em Algebra Universalis}, 81(20), 2020.
\newblock Preprint available at ArXiv:1905.12333.

\bibitem{HNBook}
Pavol Hell and Jaroslav Ne\v{s}et\v{r}il.
\newblock {\em Graphs and Homomorphisms}.
\newblock Oxford University Press, Oxford, 2004.

\bibitem{CarvalhoEgriJacksonNiven}
Catarina Carvalho, L{\'{a}}szl{\'{o}} Egri, Marcel Jackson, and Todd Niven.
\newblock On {M}altsev digraphs.
\newblock {\em Electr. J. Comb.}, 22(1):P1.47, 2015.

\bibitem{BMOOPW}
Manuel Bodirsky, Antoine Mottet, Miroslav Ol\v{s}\'ak, Jakub Opr\v{s}al,
  Michael Pinsker, and Ross Willard.
\newblock $\omega$-categorical structures avoiding height~1 identities.
\newblock {\em Transactions of the American Mathematical Society}.
\newblock Accepted.

\bibitem{JacksonKN16}
Marcel Jackson, Tomasz Kowalski, and Todd Niven.
\newblock Complexity and polymorphisms for digraph constraint problems under
  some basic constructions.
\newblock {\em Int. J. Algebra Comput.}, 26(7):1395--1433, 2016.

\end{thebibliography}


\end{document}